\documentclass[leqno]{amsart}

\numberwithin{equation}{section}

\begin{document}

\title[Hierarchy of B\"acklund Transformation Groups]
{Hierarchy of B\"acklund Transformation Groups of the Painlev\'e Systems}
\author{Masaki Suzuki, Nobuhiko Tahara and Kyoichi Takano}
\thanks{%
The third author was partially supported by Grant-in-Aid for 
Scientific Research (No. 13440054 (b), (2)), Japan Society for the Promotion 
of Science.}

\date{}

\begin{abstract}
For each Painlev\'e system $P_J$ except the first one, we have a B\"acklund 
transformation group which is a lift of an affine Weyl group. In this 
paper, we show that the B\"acklund transformation groups for 
$J=V,IV,III,II$ are successively obtained from that for $J=VI$ by the 
well known degeneration or confluence processes.
\end{abstract}

\maketitle

\section{Introduction}

The $J$-th Painlev\'e system $P_J\ (J=VI,V,IV,III,II,I)$ which is 
equivalent to the $J$-th Painlev\'e equation is the following 
Hamiltonian system
$$
 \delta_J\, q=\{H_J(q,p,t,\alpha),q\}, \quad 
 \delta_J\, p=\{H_J(q,p,t,\alpha),p\}, \leqno P_J:
$$
where $\delta_{VI}=t(t-1)d/dt$, $\delta_V=\delta_{III}=td/dt$, 
$\delta_{IV}=\delta_{II}=\delta_I=d/dt$, $\{\cdot, \cdot\}$ is a Poisson 
bracket defined by 
$$\{f,g\} =  \frac{\partial f}{\partial p} \frac{\partial g}{\partial q}
           - \frac{\partial f}{\partial q} \frac{\partial g}{\partial p},
$$ 
and the Hamiltonian $H_J=H_J(q,p,t,\alpha)$ is of the form 
\begin{eqnarray*}
H_{VI}(q,p,t,\alpha)&=&q(q-1)(q-t)p^2 -[(\alpha_0-1)q(q-1)
    +\alpha_4(q-1)(q-t) \\
   &&\quad +\alpha_3 q(q-t)]p + \alpha_2(\alpha_1+\alpha_2)(q-t) \\
   &&(\alpha_0 + \alpha_1 + 2\alpha_2 + \alpha_3 + \alpha_4 = 1),\\
H_{V}(q,p,t,\alpha)&=&q(q-1)p(p+t)-(\alpha_1+\alpha_3)qp 
    +\alpha_1 p +\alpha_2 tq \\
  &&(\alpha_0 + \alpha_1 + \alpha_2 + \alpha_3 = 1), \\
H_{IV}(q,p,t,\alpha)&=&qp(2p-q-2t) - 2\alpha_1 p - \alpha_2 q \\
  &&(\alpha_0 + \alpha_1 + \alpha_2 = 1), \\
H_{III}(q,p,t,\alpha)&=&q^2p(p-1) + q[(\alpha_0+\alpha_2)p-\alpha_0] + tp \\
  &&(\alpha_0 + 2\alpha_1 + \alpha_2 = 1), \\
H_{II}(q,p,t,\alpha)&=&\frac12p^2 - (q^2+\frac{t}{2})p - \alpha_1 q \\
  &&(\alpha_0 + \alpha_1 = 1), \\
H_I(q,p,t)&=&\frac12p^2-2q^3-tq.
\end{eqnarray*}
Notice that the Hamiltonian for $J=IV$ is slightly different from that 
in \cite{MTY} but it is of the same form as in \cite{IKSY} and \cite{T},
because we use the well known degenerations in this paper. 

The B\"acklund transformation group $W=W_J$ of Painlev\'e system 
$P_J\ (J\not= I)$ consists of birational symplectic transformations each 
of which preserves the form of the Hamiltonian $H_J$ but changes the 
parameters $\alpha=(\alpha_0,\ldots)$ as an element of an affine Weyl group. 
In other words, the elements of $W_J$ which is a lift of an affine Weyl 
group are Poisson bracket preserving differential isomorphisms of a 
differential field of functions of $q,p,\alpha$ equipped with a derivation 
defined by the system $P_J$ and $\delta_J\,\alpha_i=0$, $i=0,1,\ldots$.
Here differential isomorphism means algebraic isomorphism commuting with 
the derivation. 
The group is generated by a finite 
set of generators $s_0,s_1,\ldots$ which correspond to the simple roots of 
the affine Lie algebra (\cite{NY1}, \cite{P}). 

On the other hand, we know degenerations of Painlev\'e systems 
as the following diagram (\cite{IKSY}, \cite{M}, \cite{P}, \cite{T}):
$$
\begin{array}{ccccccccc}
{}&{}&{}&{}&P_{IV}&{}&{}&{}&{} \\
{}&{}&{}&\nearrow&{}&\searrow&{}&{}&{} \\
P_{VI}&\longrightarrow&P_V&{}&{}&{}&P_{II}
&\longrightarrow&P_I. \\
{}&{}&{}&\searrow&{}&\nearrow&{}&{}&{} \\
{}&{}&{}&{}&P_{III}&{}&{}&{}&{}
\end{array}
$$
For every $P_J\rightarrow P_K$ in the diagram, there is a change of 
parameters and variables 
\begin{eqnarray*}
&&\alpha_i=\alpha_i(A,\varepsilon)\ \ (i=0,1,\ldots), \\
&&t=t(\varepsilon,T),\quad 
q=q(A,\varepsilon,T,Q,P), \quad
p=p(A,\varepsilon,T,Q,P),
\end{eqnarray*}
between $\alpha=(\alpha_0,\alpha_1,\ldots),t,q,p$ and 
$A=(A_0,A_1,\ldots),\varepsilon, T,Q,P$. For example, in the case of 
$P_{VI}\rightarrow P_V$, 
\begin{eqnarray*}
&&\alpha_0=\varepsilon^{-1},\ \ \alpha_1=A_3,\ \ 
\alpha_2=A_2,\ \ \alpha_3=A_0-A_2-\varepsilon^{-1},\ \ \alpha_4=A_1,\\
&&t=1+\varepsilon T,\ \ (q-1)(Q-1)=1,\ \ 
(q-1)p+(Q-1)P=-A_2.
\end{eqnarray*}
Since the change of variables 
is symplectic, namely
$$
\{P,Q\}=1,\ \ \{Q,Q\}=\{P,P\}=0,
$$
the system $P_J$ is also written in the new variables $T,Q,P$ and 
parameters $A,\varepsilon$ as a Hamiltonian system 
denoted by $P_{J\to K}$. The system $P_{J\to K}$ tends to the system 
$P_K$ as $\varepsilon\to 0$ and then the process $\varepsilon\to 0$ 
in the change of parameters and variables is called a degeneration 
or confluence process from $P_J$ to $P_K$. 

In this paper, we observe how the degeneration process 
from $P_J$ to $P_K$ works on the B\"acklund transformation group 
$W_J$. The change of parameters and variables lifts the group 
$W_J$ to a group denoted again by $W_J$ each element of which is a 
differential isomorphism of a differential field of functions of 
$A=(A_0,A_1,\ldots),\varepsilon,T,Q,P$. We see that an element of 
the new $W_J$ does not converge as $\varepsilon\to 0$, in general. 
However we can verify the following theorem, which is the main assertion
of this paper.

\medskip
THEOREM. \ \ {\it For every degeneration process 
$P_J\rightarrow P_K$ except for $J=II,K=I$ in 
Painlev\'e systems, we can choose a subgroup $W_{J\to K}$ of the 
B\"acklund transformation group $W_J$ so that $W_{J\to K}$ converges to 
$W_K$ as $\varepsilon\to 0$. }

\medskip
The subgroup $W_{J\to K}$ of $W_J$ is taken 
as a group generated by reflections of $A_0,A_1,\ldots$, since 
the new parameters $A_0,A_1,\ldots$ should be the simple roots of 
an affine Weyl algebra for the system $P_K$. 

Here we notice that the same process for $P_{II} \to P_I$ can be 
followed, however we see that each generator of $W_{II\to I}$ 
converges to the identity as $\varepsilon\to 0$. 
The fact seems to suggest that the first Painlev\'e system $P_I$ has 
no nontrivial B\"acklund transformations.

Since each $W_J$ is a lift of 
an affine Weyl group corresponding to an affine Lie algebra (see next 
section), it is convenient to express the above theorem by the 
following diagram:
$$
\begin{array}{ccccccc}
{}&{}&{}&{}&W(A^{(1)}_2)&{}&{} \\
{}&{}&{}&\nearrow&{}&\searrow&{} \\
W(D^{(1)}_4)&\longrightarrow&W(A^{(1)}_3)&{}&{}&{}&W(A^{(1)}_1). \\
{}&{}&{}&\searrow&{}&\nearrow&{} \\
{}&{}&{}&{}&W(C^{(1)}_2)&{}&{}
\end{array}
$$

In Section~\ref{sec:reviewBT},
we review the B\"acklund transformation groups of 
the Painlev\'e systems $P_J\ (J\not= I)$. 
The following sections are devoted to the proof of the above theorem 
in all cases of degenerations. In these sections, we also see how 
$W_{J\to K}$ acts on the system $P_{J\to K}$. 

\section{Review of B\"acklund transformation groups}
\label{sec:reviewBT}

In this section, we give explicit forms of the generators 
$s_i$ of the B\"acklund transformation group $W$ of each Painlev\'e 
system. Each list consists of the type of affine Weyl group, Dynkin 
diagram, the fundamental relations of the generators of the group 
$W$, and the explicit forms of the generators, where the forms of 
$s_i(t)$ are omitted in the case of $s_i(t)=t$ for all $i$. 
The lists are the same as those in \cite{MTY} except the case of $J=IV$. 

\subsection{The case of $J=VI$}
\label{sec:BT6}
$$
D^{(1)}_4:\quad
\begin{picture}(60,20)(0,0)
\hbox{
\put(15,0){\circle{4}}
\put(15,10){\circle{4}}
\put(30,5){\circle{4}}
\put(45,0){\circle{4}}
\put(45,10){\circle{4}}
\put(0,10){$\alpha_0$}
\put(0,0){$\alpha_1$}
\put(26,12){$\alpha_2$}
\put(52,0){$\alpha_4$}
\put(52,10){$\alpha_3$}
\put(17,10){\line(2,-1){10}}
\put(43,10){\line(-2,-1){10}}
\put(17,0){\line(2,1){10}}
\put(43,0){\line(-2,1){10}}
}
\end{picture}
\quad (\alpha_0+\alpha_1+2\alpha_2+\alpha_3+\alpha_4=1)
$$
\smallskip
$$
W(D^{(1)}_4)=\langle s_0, s_1, s_2 ,s_3, s_4 \rangle:\ \ 
s_i^2=s_2^2=1, \quad (s_is_j)^2=1, \quad (s_is_2)^3=1.
$$
\smallskip
$$
\begin{array}{c|ccccc|cc}
& \alpha_0 &\alpha_1 & \alpha_2&\alpha_3& \alpha_4 &q & p  \\
\hline
s_0 & -\alpha_0&\alpha_1 &\alpha_2+\alpha_0&\alpha_3&\alpha_4
& q&p-\frac{\alpha_0}{q-t}\\
s_1 &\alpha_0&-\alpha_1 &\alpha_2+\alpha_1 &\alpha_3&\alpha_4 
&q &p\\
s_2 &\alpha_0+\alpha_2&\alpha_1+\alpha_2 &-\alpha_2 &\alpha_3+\alpha_2
&\alpha_4+\alpha_2&q+\frac{\alpha_2}{p} &p\\
s_3 &\alpha_0&\alpha_1 &\alpha_2+\alpha_3&-\alpha_3 &\alpha_4 
&q&p-\frac{\alpha_3}{q-1} \\
s_4 &\alpha_0&\alpha_1 &\alpha_2+\alpha_4&\alpha_3 &-\alpha_4 
&q&p-\frac{\alpha_4}{q}
\end{array}
$$

\medskip
\noindent
The last list should be read as 
\begin{align*}
  s_0(\alpha_0)=-\alpha_0,\ \  s_0(\alpha_1)=\alpha_1, \ \ &
s_0(\alpha_2)=\alpha_2+\alpha_0, \ \ s_0(\alpha_3)=\alpha_3, \ \ 
s_0(\alpha_4)=\alpha_4, \\
s_0(q)=q, \ \ &
  s_0(p)=p-\frac{\alpha_0}{q-t}
\end{align*}
and so on.

\subsection{The case of $J=V$}
\label{sec:BT5}
$$
A^{(1)}_3:\quad
\begin{picture}(60,20)(0,0)
\hbox{
\put(10,0){\circle{4}}
\put(30,10){\circle{4}}
\put(30,0){\circle{4}}
\put(50,0){\circle{4}}
\put(4,7){$\alpha_1$}
\put(26,17){$\alpha_0$}
\put(48,7){$\alpha_3$}
\put(26,-10){$\alpha_2$}
\put(15,0){\line(1,0){10}}
\put(45,0){\line(-1,0){10}}
\put(15,1){\line(2,1){13}}
\put(45,1){\line(-2,1){13}}
}
\end{picture}
\quad (\alpha_0+\alpha_1+\alpha_2+\alpha_3=1)
$$
\smallskip
$$
W(A^{(1)}_3)=\langle s_0, s_1, s_2 ,s_3 \rangle:\ \ 
s_i^2=1, \quad (s_is_{i+2})^2=1, \quad (s_is_{i+1})^3=1.
$$
\smallskip
$$
\begin{array}{c|cccc|cc}
& \alpha_0 &\alpha_1 & \alpha_2&\alpha_3& q & p  \\
\hline
s_0 & -\alpha_0&\alpha_1+\alpha_0 &\alpha_2&\alpha_3+\alpha_0 
& q+\frac{\alpha_0}{p+t}&p\\
s_1 &\alpha_0+\alpha_1 &-\alpha_1 &\alpha_2+\alpha_1&\alpha_3 
&q &p-\frac{\alpha_1}{q}\\
s_2 &\alpha_0&\alpha_1+\alpha_2 &-\alpha_2 &\alpha_3+\alpha_2
&q+\frac{\alpha_2}{p} &p\\
s_3 &\alpha_0+\alpha_3 &\alpha_1&\alpha_2+\alpha_3 &-\alpha_3 
&q&p-\frac{\alpha_3}{q-1}
\end{array}
$$

\subsection{The case of $J=IV$}
\label{sec:BT4}
$$
A^{(1)}_2:\quad
\begin{picture}(60,20)(0,0)
\hbox{
\put(10,0){\circle{4}}
\put(30,10){\circle{4}}
\put(50,0){\circle{4}}
\put(4,7){$\alpha_1$}
\put(26,17){$\alpha_0$}
\put(48,7){$\alpha_2$}
\put(15,0){\line(1,0){30}}
\put(15,1){\line(2,1){13}}
\put(45,1){\line(-2,1){13}}
}
\end{picture}
\quad (\alpha_0+\alpha_1+\alpha_2=1)
$$
\smallskip
$$
W(A^{(1)}_2)=\langle s_0, s_1, s_2 \rangle:\ \ 
s_0^2=s_1^2=s_2^2=1, \quad (s_0s_1)^3=(s_1s_2)^3=(s_2s_0)^3=1.
$$
\smallskip
$$
\begin{array}{c|ccc|cc}
& \alpha_0 &\alpha_1 & \alpha_2 & q & p  \\
\hline
s_0 & -\alpha_0&\alpha_1+\alpha_0 &\alpha_2+\alpha_0 
& q+\frac{2\alpha_0}{2p-q-2t}&p+\frac{\alpha_0}{2p-q-2t}\\
s_1 &\alpha_0+\alpha_1 &-\alpha_1 &\alpha_2+\alpha_1 
&q &p-\frac{\alpha_1}{q}\\
s_2 &\alpha_0+\alpha_2 &\alpha_1+\alpha_2 &-\alpha_2 
&q+\frac{\alpha_2}{p} &p
\end{array}
$$

\subsection{The case of $J=III$}
$$
C^{(1)}_2:\quad
\begin{picture}(60,20)(0,0)
\hbox{
\put(10,3){\circle{4}}
\put(15,1){$\Rightarrow$}
\put(30,3){\circle{4}}
\put(35,1){$\Leftarrow$}
\put(50,3){\circle{4}}
\put(4,10){$\alpha_0$}
\put(26,10){$\alpha_1$}
\put(48,10){$\alpha_2$}
}
\end{picture}
\quad (\alpha_0+2\alpha_1+\alpha_2=1)
$$
\smallskip
$$
W(C^{(1)}_2)=\langle s_0, s_1, s_2 \rangle:\ \ 
s_0^2=s_1^2=s_2^2=1, \quad (s_0s_1)^4=(s_1s_2)^4=1.
$$
\smallskip
$$
\begin{array}{c|ccc|ccc}
&\alpha_0&\alpha_1&\alpha_2&t&q&p\\
\hline
s_0&-\alpha_0&\alpha_1+\alpha_0&\alpha_2
&t&q+\frac{\alpha_0}{p}&p\\
s_1&\alpha_0+2\alpha_1&-\alpha_1&\alpha_2+2\alpha_1
&-t&q&p-\frac{2\alpha_1}{q}+\frac{t}{q^2}\\
s_2&\alpha_0&\alpha_1+\alpha_2&-\alpha_2
&t&q+\frac{\alpha_2}{p-1}&p
\end{array}
$$

\subsection{The case of $J=II$}
$$
A^{(1)}_1:\quad
\begin{picture}(50,20)(0,0)
\hbox{
\put(10,3){\circle{4}}
\put(15,1){$\Leftrightarrow$}
\put(30,3){\circle{4}}
\put(4,10){$\alpha_0$}
\put(28,10){$\alpha_1$}
}
\end{picture}
\quad (\alpha_0+\alpha_1=1)
$$
\smallskip
$$
W(A^{(1)}_1)=\langle s_0, s_1 \rangle:
\quad s_0^2=s_1^2=1.
$$
\smallskip
$$
\begin{array}{c|cc|cc}
&\alpha_0&\alpha_1&q&p\\
\hline
s_0&-\alpha_0&\alpha_1+2\alpha_0
&q+\frac{\alpha_0}{p-2q^2-t}&p+\frac{4\alpha_0q}{p-2q^2-t}
  +\frac{2\alpha_0^2}{(p-2q^2-t)^2}\\
s_1&\alpha_0+2\alpha_1&-\alpha_1&q+\frac{\alpha_1}{p}&p
\end{array}
$$

\section{Degeneration from $W_{VI}$ to $W_V$}

In this case, the degeneration process is given by 
\begin{eqnarray}
&&\alpha_0=\varepsilon^{-1},\ \ \alpha_1=A_3,\ \ 
\alpha_2=A_2,\ \ \alpha_3=A_0-A_2-\varepsilon^{-1},\ \ \alpha_4=A_1,\\
&&t=1+\varepsilon T,\ \ (q-1)(Q-1)=1,\ \ 
(q-1)p+(Q-1)P=-A_2.
\end{eqnarray}
Notice that $A_0+A_1+A_2+A_3=\alpha_0+\alpha_1+2\alpha_2+\alpha_3+\alpha_4
=1$ and the change of variables from $(q,p)$ to $(Q,P)$ is symplectic.

Each B\"acklund transformation in $W_{VI}$ given in Section~\ref{sec:BT6}
is an differential isomorphism of the differential field 
$K={\bf C}(\alpha,t,q,p)$ of rational functions of 
$\alpha=(\alpha_0,\alpha_1,\ldots,\alpha_4),t,q,p$ equipped with a 
derivation $\delta_{VI}$ defined by 
\begin{eqnarray*}
&&\delta_{VI}\,q=\{H_{VI},q\},\ \ \delta_{VI}\,p=\{H_{VI},p\},\\
&&\delta_{VI}\,t=t(t-1),\ \ \delta_{VI}\,\alpha_i=0,\ 
i=0,1,\ldots,4.
\end{eqnarray*}

Since the change of parameters and variables (3.1), (3.2) is birational, 
we can obtain the action of $W_{VI}$ on the differential field 
$K':={\bf C}(A,\varepsilon,T,Q,P)$ 
of rational functions of $A=(A_0,A_1,A_2,A_3),\varepsilon,T,Q,P$. 

Let us see the actions of the generators $s_i$, $i=0,1,2,3,4$, on 
the parameters $A_i$, $i=0,1,2,3$, and $\varepsilon$ where
$$
A_0=\alpha_0+\alpha_2+\alpha_3,\ \ A_1=\alpha_4,\ \ 
A_2=\alpha_2,\ \ A_3=\alpha_1,\ \ \varepsilon=\frac{1}{\alpha_0}. 
$$
For example, the 
action of $s_0$ is obtained as 
\begin{align*}
&s_0(A_0)=s_0(\alpha_0+\alpha_2+\alpha_3)=-\alpha_0+(\alpha_2+\alpha_0)
+\alpha_3=\alpha_2+\alpha_3=A_0-\varepsilon^{-1}, \\
&s_0(A_1)=s_0(\alpha_4)=\alpha_4=A_1,\ \ 
s_0(A_2)=s_0(\alpha_2)=\alpha_2+\alpha_0=A_2+\varepsilon^{-1}, \\
&s_0(A_3)=s_0(\alpha_1)=\alpha_1=A_3,\ \ 
s_0(\varepsilon)=s_0(1/\alpha_0)=-1/\alpha_0=-\varepsilon.
\end{align*}
Similarly we have
\begin{align*}
s_1(A_0)&=A_0+A_3,\ s_1(A_1)=A_1,\ s_1(A_2)=A_2+A_3,\ 
s_1(A_3)=-A_3,\ s_1(\varepsilon)=\varepsilon,\\
s_2(A_0)&=A_0,\ s_2(A_1)=A_1+A_2,\ s_2(A_2)=-A_2,\ 
s_2(A_3)=A_3+A_2,\\
s_2(\varepsilon)&=\frac\varepsilon{1+A_2\varepsilon},\\
s_3(A_0)&=A_2+\varepsilon^{-1},\ s_3(A_1)=A_1,\ s_3(A_2)
=A_0-\varepsilon^{-1},\ s_3(A_3)=A_3,\ s_3(\varepsilon)=\varepsilon, \\
s_4(A_0)&=A_0+A_1,\ s_4(A_1)=-A_1,\ s_4(A_2)=A_2+A_1,\ 
s_4(A_3)=A_3,\ s_3(\varepsilon)=\varepsilon. 
\end{align*}
We remark that $s_3(A_0)$ and $s_3(A_2)$ diverge as $\varepsilon\to 0$. 

Observing these relations, we take a subgroup $W_{VI\to V}$ of 
$W_{VI}$ generated by $S_0,S_1,S_2,S_3$ defined by
\begin{equation}
S_0:=s_0s_2s_3s_2s_0=s_3s_2s_0s_2s_3,\ \ 
S_1:=s_4,\ \ S_2:=s_2,\ \ S_3:=s_1.
\end{equation}
We can easily check 
\begin{eqnarray}
&&S_0(A_0)=-A_0,\ \ S_0(A_1)=A_1+A_0,\ \ 
S_0(A_2)=A_2,\\
&&S_0(A_3)=A_3+A_0,\ \ 
S_0(\varepsilon)=\frac\varepsilon{1-A_2\varepsilon},\nonumber\\
&&S_1(A_0)=A_0+A_1,\ \ S_1(A_1)=-A_1,\ \ 
S_1(A_2)=A_2+A_1, \\
&&S_1(A_3)=A_3,\ \ S_1(\varepsilon)=\varepsilon,\nonumber \\
&&S_2(A_0)=A_0,\ \ S_2(A_1)=A_1+A_2,\ \ 
S_2(A_2)=-A_2,\\
&&S_2(A_3)=A_3+A_2,\ \ 
S_2(\varepsilon)=\frac\varepsilon{1+A_2\varepsilon},\nonumber\\
&&S_3(A_0)=A_0+A_3,\ \ S_3(A_1)=A_1,\ \ 
S_3(A_2)=A_2+A_3, \\
&&S_3(A_3)=-A_3,\ \ S_3(\varepsilon)=\varepsilon,\nonumber 
\end{eqnarray}
and the generators satisfy the fundamental relations given in 
Section~\ref{sec:BT5}. In short, the group
$W_{VI\to V}=\langle S_0,S_1,S_2,S_3 \rangle$ can be 
considered to be an affine Weyl group of the affine Lie algebra of 
type $A_3^{(1)}$ with simple roots $A_0,A_1,A_2,A_3$. 

Now we investigate how the generators of $W_{VI\to V}$ act on 
$T,Q$ and $P$. We can verify 
\begin{eqnarray}
&&S_0(T)=T(1-A_0\varepsilon),\ \ S_0(Q)=Q+
\frac{A_0(1-Q(Q-1)P\varepsilon)}{P+T-T(Q-1)P\varepsilon},\\
&&S_0(P)=P\biggl(1+\frac{A_0T\varepsilon}{P+T-T(A_0+QP)\varepsilon}\biggr),
\nonumber\\
&&S_1(T)=T,\ \ S_1(Q)=Q,\ \ S_1(P)=P-\frac{A_1}{Q},\\
&&S_2(T)=T(1+A_2\varepsilon),\ \ S_2(Q)=Q+\frac{A_2}{P},\ \ 
S_2(P)=P,\\
&&S_3(T)=T,\ \ S_3(Q)=Q,\ \ S_3(P)=P-\frac{A_3}{Q-1}.
\end{eqnarray}
By comparing (3.4) -- (3.11) with the last list in Section~\ref{sec:BT5},
we see that 
our theorem holds for $W_{VI}\to W_V$. 

We notice that the system $P_{VI}$ is written in the new variables as 
$$
\delta_V\,Q=\{H_{VI\to V},Q\},\quad \delta_V\,P=\{H_{VI\to V},P\} 
\leqno P_{VI\to V}:
$$
where $\delta_V=T\partial/\partial T,\ 
H_{VI\to V}:=H_{VI}/(1+\varepsilon T)$, $H_{VI\to V}\to H_V$ as 
$\varepsilon \to 0$. We can verify that $\delta_V$ commutes with any 
element $W_{VI\to V}$, and then for any $w\in W_{VI\to V}$
$$
\delta_V\,w(Q)=\{w(H_{VI\to V}),w(Q)\},\quad 
\delta_V\,w(P)=\{w(H_{VI\to V}),w(P)\}. 
$$
\section{Degeneration from $W_V$ to $W_{IV}$}
\label{sec:deg5to4}

The degeneration in the case is given by 
\begin{eqnarray}
&&\alpha_0=A_0+\frac12\varepsilon^{-2},\ \ \alpha_1=A_1,\ \ 
\alpha_2=A_2,\ \ \alpha_3=-\frac12\varepsilon^{-2}, \\
&&t=\frac12\varepsilon^{-2}(1+2\varepsilon T),\ \ 
q=-\frac{\varepsilon Q}{1-\varepsilon Q},\\
&&p=-\varepsilon^{-1}(1-\varepsilon Q)[P-\varepsilon(A_2+QP)]. 
\nonumber 
\end{eqnarray}
Notice that $A_0+A_1+A_2=\alpha_0+\alpha_1+\alpha_2+\alpha_3=1$ and 
the transformation from $(q,p)$ to $(Q,P)$ is symplectic, however 
the change of parameters (4.1) is not one to one differently from 
the case of $P_{VI}\to P_V$. 

Since the generators of $W_{V\to IV}$ should be reflections 
of $A_0=\alpha_0+\alpha_3$, $A_1=\alpha_1$, $A_2=\alpha_2$, we choose 
them as 
\begin{equation}
S_0:=s_3s_0s_3=s_0s_3s_0,\ \ S_1:=s_1,\ \ S_2:=s_2
\end{equation}
and set $W_{V\to VI}=\langle S_0,S_1,S_2 \rangle$.
Then we immediately have 
\begin{eqnarray}
&&S_0(A_0)=-A_0,\ \ S_0(A_1)=A_1+A_0,\ \ S_0(A_2)=A_2+A_0,\\
&&S_1(A_0)=A_0+A_1,\ \ S_1(A_1)=-A_1,\ \ S_1(A_2)=A_2+A_1,\\
&&S_2(A_0)=A_0+A_2,\ \ S_2(A_1)=A_1+A_2,\ \ S_2(A_2)=-A_2.
\end{eqnarray}
However, we see that $S_i(\varepsilon)$ have ambiguities of 
signature. For example, since
$$
S_2(\varepsilon)^2=s_2(\varepsilon^2)=s_2((-1/2)/\alpha_3)
=-\frac12\frac{1}{\alpha_3+\alpha_2}
=\frac{\varepsilon^2}{1-2A_2\varepsilon^2},
$$
we can choose any one of the two branches as $S_2(\varepsilon)$. 
Among such possibilities, we take a choice as
\begin{equation}
S_0(\varepsilon)=\varepsilon(1+2A_0\varepsilon^2)^{-1/2},\ \ 
S_1(\varepsilon)=\varepsilon,\ \ 
S_2(\varepsilon)=\varepsilon(1-2A_2\varepsilon^2)^{-1/2}
\end{equation}
where $(1+2A_0\varepsilon^2)^{1/2}=1$ and 
$(1-2A_2\varepsilon^2)^{1/2}=1$ at $A_0\varepsilon^2=0$ 
and $A_2\varepsilon^2=0$ respec\-tively, or considering in the 
category of formal power series, we make a convention 
that $(1+2A_0\varepsilon^2)^{1/2}$ and $(1-2A_2\varepsilon^2)^{1/2}$ 
are formal power series of $A_0\varepsilon^2$ and $A_2\varepsilon^2$ 
with constant terms $1$ according to 
$$
(1+x)^c \sim 1+\sum_{n\geq 1}\binom{c}{n}x^n.
$$
We notice that the generators acting on parameters 
$A_0,A_1,A_2,\varepsilon$ satisfy the fundamental relations in 
Section~\ref{sec:BT4}. 

Now we observe the actions of $S_i$, $i=0,1,2$, on the variables 
$T,Q,P$. By means of (4.2), (4.7) and 
$$
S_0(t)=s_3s_0s_3(t)=t,\ \ S_1(t)=s_1(t)=t,\ \ S_2(t)=s_2(t)=t,
$$
we can easily check
\begin{eqnarray}
&&S_0(T)=(T-A_0\varepsilon)(1+2A_0\varepsilon^2)^{-1/2},\ \ 
S_1(T)=T,\\
&&S_2(T)=(T+A_2\varepsilon)(1-2A_2\varepsilon^2)^{-1/2}.
\end{eqnarray}
By (4.1), (4.2), (4.7) and the actions of $s_1,s_2$ on $q,p$, we can 
easily verify 
\begin{eqnarray}
&&S_1(Q)=Q,\quad S_1(P)=P-\frac{A_1}{Q},\\
&&S_2(Q)=Q+\frac{A_2}{P},\quad S_2(P)=P.
\end{eqnarray}
The forms of the actions $S_0=s_3s_0s_3$ on $Q$ and $P$ are complicated,
but we can see that 
\begin{equation}
S_0(Q) \to Q+\frac{2A_0}{2P-Q-2T},\ \  
S_0(P) \to P+\frac{A_0}{2P-Q-2T} 
\end{equation}
as $\varepsilon\to 0$ for arbitrarily fixed $A=(A_0,A_1,A_2),T,Q$ and $P$ 
with some generic conditions such as $2P-Q-2T\not= 0$. Here we have to 
note that, although $S_0(Q),S_0(P)$ contain formal power series of 
$A,\varepsilon$, they are analytic if $\varepsilon$ is sufficiently small 
for any fixed $A,T,Q,P$. 

By means of the above study, we define a differential field $K'$ on 
which $W_{V\to IV}=\langle S_0,S_1,S_2 \rangle$
acts as the field of rational functions 
of $T,Q,P$ whose coefficients are formal power series of 
$A_0,A_1,A_2,\varepsilon$. 
Then the action of any $w\in W_{V\to IV}$ is defined as an isomorphism 
from $K'$ to itself. 

The equations or property from (4.4) to (4.12) and
the list in Section~\ref{sec:BT4} 
show the theorem for $W_V \to W_{IV}$.

Since $\delta_V=td/dt=(1+2\varepsilon T)(2\varepsilon)^{-1}d/dT
=(1+2\varepsilon T)(2\varepsilon)^{-1}\delta_{IV}$ and 
the transformation from $(q,p)$ to $(Q,P)$ is symplectic, the 
system $P_V$ is expressed as
$$
\delta_{IV}Q=\{H_{V\to IV},Q\},\ \ \delta_{IV}P=\{H_{V\to IV},P\}
\leqno P_{V\to IV}:
$$
in the new variables,
where $H_{V\to IV}=2\varepsilon(1+2\varepsilon T)^{-1}H_V$, and
$H_{V\to IV}\to H_{IV}$ as $\varepsilon\to 0$. However $\delta_{IV}$ 
does not commutes with the elements of $W_{V\to IV}$ and then 
we have to notice that the transform of $P_{V\to IV}$ by 
$w\in W_{V\to IV}$ 
is 
\begin{eqnarray*}
&&\delta_{IV}w(Q)=\biggl\{\frac{2\varepsilon}{1+2\varepsilon T}
w\biggl(\frac{1+2\varepsilon T}{2\varepsilon}\biggr)w(H_{V\to IV}), 
w(Q)\biggr\}, \\
&&\delta_{IV}w(P)=\biggl\{\frac{2\varepsilon}{1+2\varepsilon T}
w\biggl(\frac{1+2\varepsilon T}{2\varepsilon}\biggr)w(H_{V\to IV}), 
w(P)\biggr\},
\end{eqnarray*}
which is verified by the fact that $\delta_V$ commutes with every 
$w\in W_{V\to IV}$. 

\section{Degeneration from $W_V$ to $W_{III}$}

The degeneration in this case is
\begin{eqnarray}
&&\alpha_0=A_2,\ \  \alpha_1=\varepsilon^{-1},\ \  \alpha_2=A_0,
\ \  \alpha_3=2A_1-\varepsilon^{-1},\\
&&t=-\varepsilon T,\ \ q=1+\frac{Q}{\varepsilon T},\ \ p=\varepsilon TP.
\end{eqnarray}
We see that $A_0+2A_1+A_2=\alpha_0+\alpha_1+\alpha_2+\alpha_3=1$ and 
the change of variables from $(q,p)$ to $(Q,P)$ is symplectic. 
As the case of $P_{VI} \to P_V$, the transformation given by (5.1) 
and (5.2) is birational, and we can easily obtain the actions of 
$s_i$, $i=0,1,2,3$, on the differential field 
$K'={\bf C}(A_0,A_1,A_2,\varepsilon,T,Q,P)$. 

Choose $S_i$, $i=0,1,2$, as
\begin{equation}
S_0:=s_2,\ \ S_1:=s_3s_1=s_1s_3,\ \ S_2:=s_0
\end{equation}
which are reflections of $A_0=\alpha_2$, $A_1=(\alpha_1+\alpha_3)/2$,
$A_2=\alpha_0$ respectively. 

It is easy to see that
\begin{eqnarray}
&&S_0(A_0)=-A_0,\  S_0(A_1)=A_1+A_0,\ S_0(A_2)=A_2,\ 
S_0(\varepsilon)=\frac\varepsilon{1+A_0\varepsilon},\\
&&S_1(A_0)=A_0+2A_1,\  S_1(A_1)=-A_1,\  S_1(A_2)=A_2+2A_1,\ 
S_1(\varepsilon)=-\varepsilon,\\
&&S_2(A_0)=A_0,\  S_2(A_1)=A_1+A_2,\  S_2(A_2)=-A_2,\  
S_2(\varepsilon)=\frac\varepsilon{1+A_2\varepsilon}
\end{eqnarray}
and
\begin{eqnarray}
&&S_0(T)=T(1+A_0\varepsilon),\ \ S_0(Q)=Q+\frac{A_0}{P},\ \ 
S_0(P)=P,\\
&&S_1(T)=-T,\ \ S_1(Q)=Q,\ \ 
S_1(P)=P-\frac{2A_1}{Q}+\frac{T}{Q^2}+O(\varepsilon),\\
&&S_2(T)=T(1+A_2\varepsilon),\ \ S_2(Q)=Q+\frac{A_2}{P-1},\ \ 
S_2(P)=P
\end{eqnarray}
where $O(\varepsilon)$ is a rational function of $A_i$, $i=0,1,2$,
$\varepsilon,T,Q,P$ with a factor $\varepsilon$. 
The proof of the theorem for $W_V\to W_{III}$ has thus been 
completed. 

We see that $\delta_V=td/dt=Td/dT=\delta_{III}$ and 
the system $P_V$ is written in the new variables by 
$$
\delta_{III}Q=\{H_{V\to III},Q\},\ \ 
\delta_{III}P=\{H_{V\to III},P\} \leqno P_{V\to III}:
$$
where $H_{V\to III}=H_V+QP$, which converges to $H_{III}$ as 
$\varepsilon\to 0$. Since $\delta_{III}$ commutes 
with any element of $W_{V\to III}$, the transform of $P_{V\to III}$ 
by $w\in W_{V\to III}$ is 
$$
\delta_{III}w(Q)=\{w(H_{V\to III}),w(Q)\},\ \ 
\delta_{III}w(P)=\{w(H_{V\to III}),w(P)\}.
$$

\section{Degeneration from $W_{IV}$ to $W_{II}$}
\label{sec:deg4to2}

The degeneration is 
\begin{eqnarray}
&&\alpha_0=A_0-\frac14\varepsilon^{-6},\ \ 
\alpha_1=\frac14\varepsilon^{-6},\ \ \alpha_2=A_1,\\
&&t=-\frac{1}{\sqrt 2}\,\varepsilon^{-3}(1-\varepsilon^4 T),\ \ 
q=\frac{1}{\sqrt 2}\,\varepsilon^{-3}(1+2\varepsilon^2Q),\ \ 
p=\frac{1}{\sqrt 2}\,\varepsilon P.
\end{eqnarray}
Then $A_0+A_1=\alpha_0+\alpha_1+\alpha_2=1$ and the change of variables 
from $(q,p)$ to $(Q,P)$ is symplectic. Since the change of parameters 
(6.1) is not one to one, we consider the degeneration process 
by introducing formal power series of the new parameters 
$A=(A_0,A_1),\varepsilon$. 

We choose $S_0$ and $S_1$ as 
\begin{equation}
S_0:=s_0s_1s_0=s_1s_0s_1,\ \ S_1:=s_2
\end{equation}
and put $W_{IV\to II}=\langle S_0,S_1 \rangle$. 
Note that $S_0,S_1$ are reflections of $A_0=\alpha_0+\alpha_1$,
$A_1=\alpha_2$ respectively. 

Then we can obtain
\begin{eqnarray}
&&S_0(A_0)=-A_0,\ \  S_0(A_1)=A_1+2A_0,\ \  S_0(\varepsilon)=
\varepsilon(1-4A_0\varepsilon^6)^{-1/6},\\
&&S_1(A_0)=A_0+2A_1,\ \ S_1(A_1)=-A_1,\ \  S_1(\varepsilon)=
\varepsilon(1+4A_1\varepsilon^6)^{-1/6}.
\end{eqnarray}
Here, we make the same convention as in Section~\ref{sec:deg5to4} that 
$(1-4A_0\varepsilon^6)^{-1/6}$ and 
$(1+4A_1\varepsilon^6)^{-1/6}$ respectively mean formal power series 
of $A_0\varepsilon^6$ and $A_1\varepsilon^6$ with $1$ as constant terms. 

Let $K'$ be a field of rational functions of $T,Q,P$ whose coefficients 
are formal power series of $A=(A_0,A_1),\varepsilon$. Then we can verify 
\begin{eqnarray}
&&S_0(T)\to T,\ \  S_0(Q)\to Q+\frac{A_0}{P-2Q^2-T},\\
&&S_0(P)\to P+\frac{4A_0Q}{P-2Q^2-T}+\frac{2A_0^2}{(P-2Q^2-T)^2}, 
\nonumber \\
&&S_1(T)\to T,\ \  S_1(Q)\to Q+\frac{A_1}{P},\\
&&S_1(P)\to P \nonumber
\end{eqnarray}
as $\varepsilon\to 0$. Concerning the convergence, remind the note 
in Section~\ref{sec:deg5to4}. Thus we have proved the theorem for 
$W_{IV} \to W_{II}$. 

Since $\delta_{IV}=(\sqrt{2}/\varepsilon)\delta_{II}$, the system 
$P_{IV}$ is written in the new variables as 
$$
\delta_{II}Q=\{H_{IV\to II},Q\},\ \ 
\delta_{II}P=\{H_{IV\to II},P\} \leqno P_{IV\to II}:
$$
where $H_{IV\to II}=(\varepsilon/\sqrt{2})H_{IV}$ and 
$H_{IV\to II}\to H_{II}$ as $\varepsilon\to 0$.  Notice that $\delta_{II}$ 
does not commute with elements of $W_{IV\to II}$, and the transform of 
$P_{IV\to II}$ by $w\in W_{IV\to II}$ is 
\begin{eqnarray*}
&&\delta_{II}w(Q)=\{\varepsilon w(1/\varepsilon)w(H_{IV\to II}),w(Q)\},\\
&&\delta_{II}w(P)=\{\varepsilon w(1/\varepsilon)w(H_{IV\to II}),w(P)\}.
\end{eqnarray*}

\section{Degeneration from $W_{III}$ to $W_{II}$}

In this case, the degeneration of parameters is given by 
\begin{equation}
\alpha_0=A_1,\ \ \alpha_1=\frac14\varepsilon^{-3},\ \ 
\alpha_2=A_0-\frac12\varepsilon^{-3}
\end{equation}
and that of variables is given by the composition of the following 
two transformations:
\begin{eqnarray}
&&t=-\tau^2,\ \ q=-\frac{\tau}{x},\ \ 
p=\frac{x}{\tau}(A_1+xy),\\
&&\tau=\frac{1+\varepsilon^2 T}{4\varepsilon^3},\ \ 
x=1+2\varepsilon Q,\ \ y=\frac{P}{2\varepsilon}.
\end{eqnarray}
Note that $A_0+A_1=\alpha_0+2\alpha_1+\alpha_2=1$ and the transformations 
from $(q,p)$ to $(x,y)$ and from $(x,y)$ to $(Q,P)$ are symplectic. 

Let us choose  
\begin{equation}
S_0:=(s_2s_1)^2=(s_1s_2)^2,\ \ S_1:=s_0
\end{equation}
as generators of $W_{III\to II}$. Then we see that
\begin{eqnarray}
&&S_0(A_0)=-A_0,\ \ S_0(A_1)=A_1+2A_0,\ \ S_0(\varepsilon)=-\varepsilon,\\
&&S_1(A_0)=A_0+2A_1,\ \ S_1(A_1)=-A_1,\ \ 
S_1(\varepsilon)=\varepsilon(1+4A_1\varepsilon^3)^{-1/3}.
\end{eqnarray}
In the last equation of (7.5), we have chosen $-1$ as a branch of 
$(-1)^{1/3}$ in order that $S_0^2(\varepsilon)=\varepsilon$. As in 
Sections~\ref{sec:deg5to4} and~\ref{sec:deg4to2}, we make a convention that 
$(1+4A_1\varepsilon^3)^{-1/3}$ is a formal power series of 
$A_1\varepsilon^3$ with $1$ as a constant term. 

By careful calculation, we can verify 
\begin{eqnarray}
&&S_0(T)=T,\ \ S_0(Q)\to Q+\frac{A_0}{P-2Q^2-T},\\
&&S_0(P)\to P+\frac{4A_0Q}{P-2Q^2-T}+\frac{2A_0^2}{(P-2Q^2-T)^2},
\nonumber \\
&&S_1(T)\to T,\ \ S_1(Q)\to Q+\frac{A_1}{P},\ \ 
S_1(P)\to P
\end{eqnarray}
as $\varepsilon\to 0$ for arbitrarily fixed $A_0,A_1,T,Q,P$. 
Thus we have proved the theorem for $W_{III}\to W_{II}$. 

We see that $\delta_{III}=(1+\varepsilon^2 T)(2\varepsilon^2)^{-1}
\delta_{II}$ and the system 
$P_{III}$ is written in the new variables as 
$$
\delta_{II}Q=\{H_{III\to II},Q\},\ \ 
\delta_{II}P=\{H_{III\to II},P\} \leqno P_{III\to II}:
$$
where $H_{III\to II}=(2\varepsilon^2)(1+\varepsilon^2 T)^{-1}H_{III}$ 
and $H_{III\to II}\to H_{II}$ as $\varepsilon\to 0$.
We notice that  $\delta_{II}$ 
does not commute with elements of $W_{III\to II}$, and the transform of 
$P_{III\to II}$ by $w\in W_{III\to II}$ is 
\begin{eqnarray*}
&&\delta_{II}w(Q)=\biggl\{\frac{2\varepsilon^2}{1+\varepsilon^2 T}
w\biggl(\frac{1+\varepsilon^2 T}{2\varepsilon^2}\biggr)w(H_{III\to II}), 
w(Q)\biggr\}, \\
&&\delta_{II}w(P)=\biggl\{\frac{2\varepsilon^2}{1+\varepsilon^2 T}
w\biggl(\frac{1+\varepsilon^2 T}{2\varepsilon^2}\biggr)w(H_{III\to II}), 
w(P)\biggr\}.
\end{eqnarray*}

\bigskip
{\bf Acknowledgement}.  The problem studied in this paper was proposed 
by Prof. Kenji Iohara in Kobe University. The authors thank 
him not only for the proposal but also for helpful discussions.

\bigskip
\leftline{Masaki SUZUKI}
\smallskip
\leftline{Graduate School of Science and Technology}
\leftline{Kobe University}
\leftline{Rokko, Kobe 657-8501}
\leftline{Japan}
\leftline{E-mail: suzuki@math.kobe-u.ac.jp}
\medskip
\leftline{Nobuhiko TAHARA}
\smallskip
\leftline{Graduate School of Science and Technology}
\leftline{Kobe University}
\leftline{Rokko, Kobe 657-8501}
\leftline{Japan}
\leftline{E-mail: tahara@math.kobe-u.ac.jp}
\medskip
\leftline{Kyoichi TAKANO}
\smallskip
\leftline{Department of Mathematics}
\leftline{Faculty of Science}
\leftline{Kobe University}
\leftline{Rokko, Kobe 657-8501}
\leftline{Japan}
\leftline{E-mail: takano@math.kobe-u.ac.jp}


\begin{thebibliography}{9}
\bibitem{IKSY}
K. Iwasaki, H. Kimura, S. Shimomura and M. Yoshida, From 
Gauss to Painlev\'e, Vieweg, 1991. 
\bibitem{M}
T. Masuda, On a class of algebraic solutions to the Painlev\'e 
VI equation, its determinant formula and coalescence cascade, Funkcial. 
Ekvac., {\bf 46} (2003), 121--171.
\bibitem{N}
M. Noumi, Painlev\'e equations, Asakura Shoten, Tokyo, 2000 
(in Japanese). 
\bibitem{MTY}
M. Noumi, K. Takano and Y. Yamada, B\"acklund transformations 
and the manifolds of Painlev\'e systems, Funkcial. Ekvac., {\bf 45}
(2002), 237--258.
\bibitem{NY1}
M. Noumi and Y. Yamada, Affine Weyl groups, discrete dynamical 
systems and Painlev\'e equations, Commun. Math. Phys., {\bf 199} (1998), 
281--295. 
\bibitem{NY2}
M. Noumi and Y. Yamada, Affine Weyl group symmetries in Painlev\'e 
type equations, in Toward the exact WKB analysis of differential 
equations, linear or non-linear (Eds. C.J. Howl, T. Kawai, Y. Takei),
245--259, Kyoto Univ. Press, Kyoto, 2000. 
\bibitem{O}
K. Okamoto, Studies on the Painlev\'e equations, I, Ann. Mat. Pura 
Appl., {\bf 146} (1987), 337--381; II, Jap. J. Math., {\bf 13} (1987), 47--76; 
III, Math. Ann., {\bf 275} (1986), 221--256; IV, Funkcial. Ekvac., {\bf 30}
(1987), 305--332.
\bibitem{P}
P. Painlev\'e, Sur les \'equations diff\'erentielles du second 
ordre \`a points critiques fixes, C. R. Acad. Sci. Paris, {\bf 143} (1906),
1111--1117.
\bibitem{T}
K. Takano, Confluence processes in defining manifolds for 
Painlev\'e systems, Tohoku Math. J., {\bf 53} (2001), 319--335.
\end{thebibliography}
\end{document}